\documentclass[11pt, a4paper]{amsart}

\usepackage[a4paper,
            left=2cm,
            right=2cm,
            top=3cm,
            bottom=3cm]{geometry}
\usepackage{color}
\newtheorem{teo}{Theorem}[section]
\newtheorem{coro}[teo]{Corollary}
\newtheorem{lem}[teo]{Lemma}
\newtheorem{pro}[teo]{Proposition}

\newtheorem{defn}[teo]{Definition}

\newtheorem{oss}[teo]{Remark}

\newcommand{\D}{\mathbb{D}}
\newcommand{\B}{\mathbb{B}}
\newcommand{\HH}{\mathbb{H}}
\newcommand{\N}{\mathbb{N}}
\newcommand{\C}{\mathbb{C}}
\newcommand{\rr}{\mathbb{R}}
\newcommand{\s}{\mathbb{S}}

\DeclareMathOperator{\Hol}{Hol}

\newcommand{\de}{\partial_c}
\newcommand{\p}{\partial}

\DeclareMathOperator{\ext}{ext}
\DeclareMathOperator{\supp}{supp}


\title{Quaternionic Carleson measures}

\author{Nikolaos Chalmoukis}
\address{Dipartimento di Matematica e Applicazioni, Universit\'a degli studi di Milano Bicocca, via Roberto Cozzi, 55 20125, Milano, Italy}
\email{nikolaos.chalmoukis@unimib.it}

\author{Giulia Sarfatti}
\address{Dipartimento di Ingegneria Industriale e Scienze Matematiche, Universit\`a Politecnica delle Marche,               
	Via Brecce Bianche 12,  60131, Ancona, Italy}
	\email{g.sarfatti@univpm.it}

	\thanks{The first author is a member of the INdAM group GNAMPA and is partially supported by the grant INdAM-GNAMPA Project, CUP E53C25002010001, "Transferring Harmonic Analysis between Discrete Structures and Manifolds" and by the research grant "Yields of the ubiquity and the geometry of inner
functions (YoungInFun)", PID2024-160326NA-I00. The second author is partially supported by 
GNSAGA-INdAM via the project ``Hypercomplex function theory and
applications'', by the MUR PRIN 2022 ``Interactions between
Geometric Structures and Function Theories'' and by Finanziamento Premiale FOE 2014 “Splines for accUrate NumeRics: adaptIve
models for Simulation Environments”.
}


\begin{document}

\begin{abstract}

In this paper we provide a general construction of a quaternionic Banach space of slice regular functions from a given Banach space of holomorphic functions, which we call its quaternionic lift. To the best of our knowledge, this construction encompasses all known examples of quaternionic Banach spaces of slice regular functions in the literature. Our main result is a characterization of Carleson and vanishing Carleson measures for such quaternionic Banach function spaces in terms of the corresponding Carleson measures of the underlying holomorphic function space. This offers a unified approach to a problem that so far has been treated on a case-by-case basis. 
\end{abstract}
\maketitle

\section{Introduction}

Carleson measures are a central object in complex analysis, introduced by Carleson in \cite{Carleson} to solve the well-known Corona Problem. For a Banach space of holomorphic functions $X $  on the complex unit disc, Carleson measures can be defined as the measures $\mu$ for which $X $   embeds continuously into $L^p(\D, \mu)$ for a given $p$. They have applications in several areas of complex analysis and operator theory. For instance, they can be used to characterize the continuity of Toeplitz, Hankel and composition operators \cite{Zhu}, and they appear in the study of interpolation problems.
When, in particular, the embedding of $X $  in $L^p(\D, \mu)$ is compact one obtains the so-called {\em vanishing Carleson measures}, which, among other things, are useful for determining when the previously mentioned operators are compact. 

Slice regularity is a notion of holomorphicity in the quaternionic setting, introduced in \cite{GSAdvances}, which gave rise to a rich theory, still in rapid development (see \cite{libroGSS} and references therein for a detailed account of the theory). For the purpose of the present paper, we can identify the class $\mathcal{SR}(\B)$ of slice regular functions on the quaternionic unit ball $\B$ with power series $q \mapsto \sum_{n \ge 0}q^na_n$ with quaternionic coefficients which converge uniformly on compact subsets of $\B$. 
Function spaces in this setting have been extensively studied from various perspectives by several authors in recent years (see, for instance, \cite{LibroQH} and references therein). 
In this context, Carleson measures were introduced in \cite{Hankel} for the Hardy space $H^2(\B)$. In \cite{SabadiniSaracco} the authors study Carleson measures for the Hardy spaces $H^p(\B)$ and the Bergman spaces $\mathcal A^p(\B)$. In \cite{Kumar}, Carleson measures for some quaternionic Dirichlet type spaces are characterized. The techniques appearing in those articles are tailored to each space separately, hence they cannot be immediately generalized to other spaces. 

In this paper we introduce a method which allow as to  give a general result to characterize Carleson and vanishing-Carleson measures for a wide class of quaternionic Banach spaces of slice regular functions constructed as a {\em quaternionic lift} of a complex Banach space (Definition \ref{quatliftintro}). 
We point out that all the quaternionic Banach spaces of slice regular functions that we are aware of can be obtained through this procedure, and therefore our main theorem  simultaneously  generalizes all results obtained so far for Carleson measures in the quaternionic setting.

We shall now describe in more detail the content of the paper. Let $X_\D$ be a complex Banach space of analytic functions on the complex unit disc $\D$ endowed with the norm $\|\cdot \|_{X_\D}$. 
In what follows, we assume that $X_\D$ embeds continuously into the space of holomorphic functions on the complex unit disc $\D$ with respect to the locally uniform topology. 

\noindent Let $f:\B \to \mathbb{H}$ be a slice regular function and let $\s$ denote the two-dimensional sphere of imaginary units in the quaternions. According to the so-called Splitting Lemma (see Lemma \ref{split} below), for any $I$ in $\s$, 
and for any $J\in \s$ orthogonal to $I$ the restriction of $f$ to $\B_I:=\B \cap (\rr+\rr I)$ can be written as $f(x+yI)=F(x+yI)+G(x+yI)J$, where we can identify $\B_I$ with $\D$, and $F, G$ with holomorphic functions on $\D$, via the natural isomorphism $x+yI \mapsto x+yi$. \\
We set $X_{\B_I}$ to be the collection of slice regular functions $f$ such that both $F$ and $G$ belong to $X_\D$, and we define
\begin{equation}\label{normanuovaintro}
	\|f\|_{X_{\B_I}}:=\sqrt{\|F\|^2_{X_\D}+\|G\|^2_{X_\D}}. 
\end{equation}

\noindent We point out that the previous definition does not depend on the choice of the imaginary unit $J$ orthogonal to $I$. Different choices of $J $ produce functions $G $ that differ by multiplication by a unimodular constant in $\mathbb{R}+ \mathbb{R} I $.

\begin{defn}[Quaternionic lift]\label{quatliftintro}
	Let $X_\D$ be a complex Banach space of holomorphic functions on $\D$ endowed with a norm  $\| \cdot \|_{X_\D}$.
We define the {\em quaternionic lift} of $X_\D$ as 
\[X_\B:=\{ f \in \mathcal{SR}(\B) \ : \ f\in X_{\B_I} \ \text{for all $I\in \s$, and } \   \|f\|_{X_\B}:=\sup_{I\in \s}\|f\|_{X_{\B_I}}<+\infty\}.\]
\end{defn}

\noindent It is not difficult to see that $\| \cdot \|_{X_\B}$ is a norm. Its completeness follows from Lemma \ref{lemma} below. Thus $X_\B$ is a quaternionic Banach space.
  
\noindent The only hypothesis on the complex space $X_\D$ that we need to prove our characterization results is that the operator $\mathcal J: X_\D \to X_\D$, defined as
\[\mathcal J F (z):=\overline{F(\overline z)} \quad \text{for any $F\in X_\D$},\]
is bounded.

\noindent As a special instance of this construction, we may take $X_\D$ to be the complex Hilbert space of holomorphic functions $F:\D \to \C$, $F(z)=\sum_{n=0}^{+\infty}z^n\alpha_n$, where the monomials form an orthogonal complete system, endowed with the norm 
\begin{equation}\label{dirC}
	\|F\|^2_{X_{\D}}:= \sum_{n=0}^{+\infty} c_n \left |\alpha_n\right|^2<+\infty,
\end{equation}
where $c_n \in \rr^+$ for any $n\in \N$ are such that $\liminf_{n}\sqrt[n]{c_n}\ge 1.$ Observe that in this case the operator $\mathcal J$ is unitary.

\noindent Consider now the function $f(q)=\sum_{n=0}^{+\infty}q^na_n$, slice regular on $\B$. If $a_n\in \HH $ decomposes, with respect to  the orthogonal imaginary units $I$ and $J$, as $a_n=\alpha_n+\beta_n J$ with $\alpha_n,\beta_n\in \rr+\rr I$  for every $n \in \N$, then the restriction of $f$ to $\B_I$ splits with respect to $I,J\in \s$, as 
\[f(z)=\sum_{n=0}^{+\infty}z^n\alpha_n+\sum_{n=0}^{+\infty}z^n\beta_n \ J=:F(z)+G(z)J \]
for any $z\in \B_I$.
Hence, by definition \[\|f\|_{X_{\B_I}}=\sqrt{\|F\|^2_{X_{\D}}+\|G\|^2_{X_{\D}}},\]
so that
\[ \|f\|_{X_{\B_I}}=\sqrt{\sum_{n=0}^{+\infty} c_n \left |\alpha_n\right|^2+\sum_{n=0}^{+\infty} c_n \left |\beta_n\right|^2}= \sqrt{\sum_{n=0}^{+\infty} c_n  \left |a_n\right|^2},\]
which does not depend on $I$. Therefore
\[\|f\|_{X_{\B}}=\sup_{I\in\s}\|f\|_{X_{\B_I}}=\sqrt{\sum_{n=0}^{+\infty} c_n  \left |a_n\right|^2}.\]

\noindent Within the above correspondence lie, for instance, the cases of the quaternionic Hardy space $H^2(\B)$, Dirichlet space $\mathcal D(\B)$, and Bergman space $\mathcal A^2(\B)$. See, e.g., \cite{LibroQH} for the definition of such spaces.

In view of the geometric properties of slice regular functions, starting from a measure $\mu$ on the quaternionic unit ball $\B$, a natural way to induce a measure on the complex unit disc $\D$ is by \textquotedblleft projecting" all the mass of $\mu$ on a fixed slice (see \cite{Hankel}). 
\begin{defn}
	Let $\mu$ be a measure on $\B$ and let us identify the complex unit disc $\D$ with the {slice $\B_i = \mathbb B \cap (\rr+ \rr i)$. } 
	We define the measure $\mu^s$ on $\D$ as 
	\[\mu^s(E)=\mu (E\cap \mathbb{R})+\mu(\{x+yI : \ y>0, \ I\in \s, \  x+yi \in E\}),\]
	for any $E \subseteq \D$.
	
\end{defn}  

We can now state our main result,  characterizing Carleson measures in the quaternionic setting. 

\begin{teo}\label{mainintro}
	Let $p\ge 1$.
	If the operator $\mathcal J$ is bounded on the complex Banach space $X_\D$, then
	a finite nonnegative Borel measure $\mu$ is Carleson for the quaternionic Banach space $X_\B$, if and only if $\mu^s$ is Carleson for $X_\D$,
	i.e. 
	\[X_\B \subseteq L^p(\B, d\mu) \quad \text{ continuously} \quad \iff \quad  X_\D \subseteq L^p(\D, d\mu^s) \quad \text{ continuously}.\]
	Furthermore,
	$\mu$ is vanishing Carleson for $X_\B$, if and only if $\mu^s$ is vanishing Carleson for $X_\D$,
	i.e. 
	\[X_\B \subseteq L^p(\B, d\mu) \quad \text{compactly} \quad \iff \quad  X_\D \subseteq L^p(\D, d\mu^s) \quad \text{compactly} .\]
\end{teo}

When we want to stress the choice of the particular space $L^p$ in which the considered Banach space is contained, we will use the term {\em $p$-Carleson measure}.

As an application we derive explicit characterizations of continuous and compact embeddings for some relevant spaces of quaternionic functions. 
\subsection{Hardy spaces $H^p(\B)$}
The quaternionic Hardy spaces were introduced in \cite{hardy} as the spaces of slice regular functions on $\B$, bounded with respect to the norm
\begin{equation*}
	\|f\|_{H^p(\B)}:=\sup_{I\in \s}\|f\|_{H^p(\B_I)},\end{equation*}
where
\begin{equation}\label{normavecchia}
	\|f\|^p_{H^p(\B_I)}:=\lim_{r\to 1^-}\frac{1}{2\pi}\int_0^{2\pi}|f(re^{It})|^pdt.
\end{equation}
Since, for any $a\in \rr^+$, the quantities
\[ (|t|^2+|s|^2)^{\frac{1}{2}}\quad \text{and} \quad (|t|^a+|s|^a)^{\frac{1}{a}} \]
are comparable with a constant depending on $a$ only, we get that the norm \eqref{normavecchia} is equivalent to the norm \eqref{normanuovaintro} and therefore $H^p(\B)$ coincides with the quaternionic lift 
of the classical complex Hardy space $H^p(\D)$ (see Definition \ref{quatliftintro}). 

For this class of function spaces the operator $\mathcal J$ is easily seen to be an isometry, hence we can apply Theorem \ref{mainintro} obtaining the characterization of $p$-Carleson and vanishing $p$-Carleson measures for $H^p(\B)$.

For any $q=re^{Jt}\in \B$, with $r\ge 0,t\in \rr,$ and $J\in \s$, we denote by 
\[S(q)=\{\rho e^{I\alpha} \in \B : |\alpha -t|\le 1-r, \ 0 < 1-\rho \le 1 -r,\ I \in \s \},\] 
the {\em symmetric box} in $\B$, indexed by $q$, which is independent of the particular $J$. 

\begin{coro}\label{hp}
Let $1\le p <+\infty$. A finite positive Borel measure $\mu$ on $\B$ is:
	\begin{itemize}
		\item[(i)] $p$-Carleson for $H^p(\B)$ if and only if for every $ q \in \B$, the measure of the symmetric box $S(q)$ satisfies $\mu(S(q)) \le c(1-|q|)$;
		\item[(ii)] vanishing $p$-Carleson for $H^p(\B)$ if and only if for every $ q \in \B$, the measure of the symmetric box $S(q)$ satisfies $\lim_{|q|\to1^-}\frac{\mu(S(q))}{1-|q|}=0$.
	\end{itemize} 
\end{coro}	  
\begin{proof}
If $q=x+yI$, with $y\ge 0$, denote by $z=x+yi$ and by $S_i(z)$ the Carleson box in the complex unit disc
\[S_i(z)=\{\rho e^{i\alpha} \in \B_i : |\alpha -t|\le 1-r, \ 0 < 1-\rho \le 1 -r\}.\] 
By definition of the measure $\mu^s$  
\[\mu(S(q))=\mu^s(S_i(z)).\]
Then the corollary follows by the corresponding results in the complex case, which can be found, for example, in Theorem 9.3 in \cite{Duren} and Theorem 8.2.5 in \cite{Zhu}.	
\end{proof}

We point out that part $(i)$ of Corollary \ref{hp} has already been proved in \cite{Hankel} for the Hilbert case and in \cite{SabadiniSaracco} for $p \neq 2$.

\subsection{Weighted Besov spaces $B^p_\alpha$}	
Let $\alpha >-1$, $p\ge \alpha + 1$. The weighted Besov space $B_\alpha^p(\D)$ is the space of holomorphic functions on the unit disc such that 
\[\|F\|^p_{B_\alpha^p(\D)}:=|F(0)|^p+\int_{\D}|F'(z)|^p(1-|z|^2)^{\alpha}dA(z)<+\infty,\]
where $dA(z)$ is the standard area element in $\D$.
We will denote by $B^p_\alpha(\B)$ the quaternionic lift of $B_\alpha^p(\D)$.

If the splitting of a slice regular function $f$ is $f=F+GJ$ on the slice $\B_I$, then its {\em slice derivative} $\partial_c f$ (see Section \ref{prelim}) splits on $\B_I$ as
\[\partial_c f(x+yI)=F'(x+yI)+G'(x+yI)J, \quad \text{and} \quad |\partial_c f(x+yI)|^2=|F'(x+yI)|^2+|G'(x+yI)|^2.\]
Hence it is not difficult to see that $B^p_\alpha(\B)$ coincides with the space defined in \cite[Definition 4.5]{Besov}.

Notice that the operator $\mathcal J$ is an isometry  for $B_\alpha^p(\D)$, and hence we can apply Theorem \ref{mainintro} obtaining the characterization of $p$-Carleson and vanishing $p$-Carleson measures for $B^p_\alpha(\B)$ from their complex counterparts.
To state the result, first, we need to define symmetric boxes constructed from an open subset of the boundary $\partial \D$ of $\D$.

\begin{defn}
	Let $O$ be an open subset of $\partial \D$ identified with $\partial \B_{i}$. If $O=\cup_{n \in \mathbb N}A_n$ with $A_n$ open disjoint arcs in $\partial\B_{i}$, we set $q_n\in \mathbb \B_{i}$ to be such that $\partial S(q_n) \cap \partial \B_{i}$ coincides with the closure of ${A_n}$ for any $n \in \mathbb N$. We define the symmetric box on the open set $O$ as
	\[S(O)=\cup_{n\in \mathbb N}S(q_n).\]
\end{defn}

\noindent We next recall the definition of classical (non-analytic) Besov spaces.
\begin{defn}
	Let $p>0$ and $0<\sigma<1$. The classical (non-analytic) Besov space $\mathcal B^p_\sigma$ consists
	of all real functions $f$ in $L^p(\partial \D, d\theta)$ such that	
	\[\|f\|^p_{\mathcal B^p_\sigma}:= \|f\|^p_{L^p(\partial \D)}+\int_{0}^{2\pi}\int_{0}^{2\pi}\frac{|f(e^{i\theta})-f(e^{i(\theta+\varphi)})|^p}{|1-e^{i\varphi}|^{1+\sigma p}}d\theta d\varphi<+\infty.\]
\end{defn}
\noindent Associated with the Besov spaces there is a notion of capacity. 
\begin{defn}
	Let $O$ be an open subset of $\partial \D$. The Besov $\mathcal{B}^p_\sigma $  capacity of $O $  is defined as
	\[{\rm cap}(O; \mathcal B^p_\sigma ):=\inf\{\|f\|^p_{\mathcal B^p_\sigma} \ :  \ f\ge 1 \ \ \text{on $O$}\}.\]
\end{defn}
\noindent We can now state the characterization of $p$-Carleson measures for the quaternionic weighted Besov spaces. 
\begin{coro}
	Let $\alpha>-1$ and $p\ge \alpha+1$. A finite positive Borel measure $\mu$ on $\mathbb B$ is $p$-Carleson for $B^p_\alpha(\B)$ if and only if there exists a positive constant $C$ such that:
	\begin{itemize}
		\item[(i)] for $\max\{1,\alpha+1\}<p\le\alpha+2$,
		\[\mu(S(O))\le C \ {\rm cap} (O,\mathcal B ^p_{1-(\alpha+1)/p}) \quad \text{for any open set $O\subseteq \partial \mathbb{D}$;} \]
		\item[(ii)] for $1<p=\alpha + 1 \le 2$,
		\[\mu (S(q))\le C(1-|q|) \quad \text{for any $q \in \B$.}\]
	\end{itemize}
	Furthermore, for $p>\alpha+2 $ $\mu $ is always a $p$-Carleson measure. 
\end{coro}
\begin{proof}
	The proof follows from Theorem 1 in \cite{Wu}, taking into account that $\mu(S(O))=\mu^s(S(O)\cap \mathbb B_{i})$ and that, in the complex unit disc, the Carleson condition can be given equivalently in terms of Carleson {\em tents} or in terms of Carleson boxes (see, e.g., \cite{Nicola}). 
\end{proof}

Concerning vanishing $p$-Carleson measures, one can use in a similar way results in the complex setting to obtain a characterization for the quaternionic weighted Besov spaces. We give here, as a particular example, the characterization for the quaternionic Dirichlet space, which corresponds to $B^2_0(\B)$.
\begin{coro}
	A finite positive Borel measure $\mu$ on $\mathbb B$ is vanishing 2-Carleson for $B^2_0(\B)$ if and only if 
		\[\sup \frac{\mu(S(O))}{{\rm cap(O,\mathcal B^2_{1/2})}}\to 0, \quad \text{as} \quad t\to 0\]
	where the supremum is taken over all open subsets  $O \subseteq \partial \D$ such that for any $q\in S(O)$ $|q|>1-t$.
\end{coro} 

\begin{proof}
	The proof follows, as in the previous cases, from Theorem 5.3.7 in \cite{Toth}.
\end{proof}
It is clear that with the same procedure one can obtain characterizations of Carleson and vanishing Carleson measures for a wide class of spaces. In particular, we recover the characterization of Carleson measures for the quaternionic Bergman spaces as in \cite{SabadiniSaracco}, as well as the results in \cite{Kumar}.	
\begin{oss}
	We point out that we could not find in the literature a characterization analogous to the one in Theorem 5.3.7 in \cite{Toth} of $p$-vanishing Carleson measures for $B^p_\alpha(\D)$ for general $p,\alpha$ such that $\max\{1,\alpha+1\}<p\le\alpha+2$.
\end{oss}

 The paper is structured as follows: in Section \ref{prelim} we recall the definition of slice regularity and some basic properties of slice regular functions needed in the sequel. Moreover we prove some preliminary results concerning the space $X_\B$ constructed as the quaternionic lift from a complex Banach space $X_\D$. The rest of the paper is dedicated to prove Theorem \ref{mainintro}: in Section \ref{sec3} we characterize $p$-Carleson measures (see Theorem \ref{main}) and in Section \ref{sec4} we characterize vanishing $p$-Carleson measures (see Theorem \ref{main2})).
\section{Preliminary results}\label{prelim}

Let us begin by recalling the definition of slice regular functions over the quaternions $\HH=\rr+i\rr+j\rr+k\rr$, together with some basic properties that will be used in the sequel. We will restrict our attention to functions defined on the quaternionic unit ball $\B=\{q\in \HH \ : \ |q|<1\}$.
We refer to the book \cite{libroGSS} for all details and proofs.

Let $\s$ denote the two-dimensional sphere of imaginary units of $\HH$, $\s=\{q\in \HH \, | \, q^2=-1\}$. Then
one can ``slice'' the space $\HH$ in copies of the complex plane that intersect along the real axis,  
\[ \HH=\bigcup_{I\in \s}(\rr+\rr I),  \hskip 1 cm \rr=\bigcap_{I\in \s}(\rr+\rr I),\]
where $\C_I:=\rr+\rr I\cong \C$, for any $I\in\s$.  
Each element $q\in \HH$  can be expressed as $q=x+yI_q$, where $x,y$ are real (if $q\in\rr$, then $y=0$) and $I_q$ is an imaginary unit.  
To have a unique decomposition (outside the real axis) we choose $y\geq 0$.
 
A function $f:\B \to \HH$ is called {\em slice regular} if for any $I\in\s$ the restriction $f_I$ of $f$ to $\B_I=\B\cap \mathbb C_I$
is {\em holomorphic}, i.e. it has continuous partial derivatives and it is such that
\[\overline{\p}_If_I(x+yI):=\frac{1}{2}\left(\frac{\p}{\p x}+I\frac{\p}{\p y}\right)f_I(x+yI)=0\]
for all $x+yI\in \B_I$.

It can be proved that
a function $f$ is slice regular on $\B$
if and only if $f$ has a power series expansion
$f(q)=\sum_{n\ge 0}q^na_n$ with quaternionic coefficients converging uniformly on compact subsets in $\B$.

The {\em slice (or Cullen) derivative} of a function $f$ which is slice regular on $\B$, is the slice regular function defined as 
\[\de f(x+yI)=\frac{1}{2}\left(\frac{\p}{\p x} -I\frac{\p}{\p y}\right)f_I(x+yI).\]

\noindent An important property of slice regular functions on $\B$ is the following.
\begin{teo}[Representation and Extension Formula]\label{RF}
Let $f:\B\to \HH$ be a regular function and let $x+y\s\subset \B$. Then, for any $I,J\in\s$,
\[f(x+yJ)=\frac{1}{2}[f(x+yI)+f(x-yI)]+J \frac I 2 [f(x-yI)-f(x+yI)].\]
	Moreover, the previous formula allows one to extend any holomorphic function $F_I:\B_I\to \HH$ to a unique slice regular function $\ext(F_I):\B\to \HH,$
	\[\ext(F_I)(x+yJ)= \frac{1}{2}[F_I(x+yI)+F_I(x-yI)]+J \frac I 2 [F_I(x-yI)-F_I(x+yI)].\] 
\end{teo}

A basic result that establishes a relation between slice regular functions and holomorphic functions of one complex variable is the following.
\begin{lem}[Splitting Lemma]\label{split}
Let $f$ be a slice regular function on $\B$. Then for any $I\in\s$ and for any $J\in \s$, $J\perp I$ there exist two holomorphic functions $F,G:\B_I=\B\cap \C_I\to \C_I$ such that
\[f(x+yI)=F(x+yI)+G(x+yI)J\]
for any $x+yI\in\B_I$.
\end{lem} 

Given a slice regular function on $\B$ $f(q)=\sum_{n \ge 0}q^na_n$, its {\em regular conjugate} is the slice regular function defined as 
\[f^c(q)=\sum_{n \ge 0}q^n\overline{a_n}.\]
\noindent By direct computation, it is not difficult to prove the following expression for the splitting components of the regular conjugate.
\begin{pro}\label{splitconj}
Let $f:\B\to \HH$ be a slice regular function and let $I,J \in \s$ be two orthogonal imaginary units. If $f$ splits on $\B_I$ with respect to $J$, as $f=F+GJ$, then the splitting of $f^c$ on $\B_I$ with respect to $J$ is given by
\[f^c(x+yI)=\overline{ F(x-yI)}-G(x+yI)J.\]
\end{pro} 

In what follows $X_\D$ will denote a Banach space of holomorphic functions on the complex unit disc, equipped with a norm $ \Vert \cdot \Vert_{X_\mathbb{D}} $ and
$X_\B$ will be its {\em quaternionic lift} constructed as in Definition \ref{quatliftintro}.
Let us prove some preliminary results on the space $X_\B$.

\begin{lem}\label{uoc}
	$X_\B$ embeds continuously into $\mathcal{SR}(\B)$ with respect to the topology of uniform convergence on compact sets.
\end{lem}
\begin{proof}
Let $f\in X_\B$ and let $r<1$. Then, thanks to the Representation Formula \ref{RF}, 
	\begin{equation*}
		\begin{aligned}
			\sup_{q \in r\B}|f(q)|&=\sup_{x+yJ \in r\B} \left|\frac{1-Ji}{2}f(x+yi)+\frac{1+Ji}{2}f(x-yi)\right|
\le 2\sup_{z\in r\B_i}|f(z)|\\
			&=2\sup_{z\in r\B_i}\sqrt{|F(z)|^2+|G(z)|^2}\le 2 \left(\sup_{z\in r\B_i}|F(z)|+\sup_{z\in r\B_i}|G(z)|\right)
	\end{aligned}
	\end{equation*}
	where $F,G$ are the splitting components of $f$ on the slice $\B_i$ with respect to $j$. Identifying $\B_i$ with $\D$, we have that $F,G \in X_\D$, which embeds continuously in $\Hol(\D)$. Therefore, for any $r<1$ there exists a constant $c_r$ depending only on $r$  such that 
	\[\sup_{z\in r\D} |F(z)| \le c_r \|F\|_{X_\D}, \quad \text{and} \quad \sup_{z\in r\D} |G(z)| \le c_r \|G\|_{X_\D},\]
	which implies that 
	\[\sup_{q\in r\B}|f(q)| \le 2c_r(\|F\|_{X_\D}+\|G\|_{X_\D}) \le 4 c_r \|f\|_{X_\B},\]
	that is $X_\B$ embeds continuously into $\mathcal{SR}(\B)$ with respect to the uniform convergence on compact sets.
\end{proof}

Let us investigate the relation between the norms $\| \cdot \|_{X_{\B_I}}$ on different slices and the norm on $X_\B$.
\begin{oss}
If  $I\neq J$, in general, the two norms $\|\cdot \|_{X_{\B_I}}$ and $\|\cdot \|_{X_{\B_J}}$ do not coincide. 
As an easy example, one can consider the norms $\|\cdot\|_{H^p(\B_i)}$ and $\|\cdot\|_{H^p(\B_j)}$ of the function $f(q)=q+i$ for $p\neq 2$. 
\end{oss}

\noindent 
A function $f\in \mathcal{SR}(\B)$ such that $f(\B_I)\subseteq \C_I$ for any $I\in \s$ is said to be {\em slice preserving}.
For slice preserving function, the norms on different slices are actually the same. In the general case norms on different slices can be compared.

\begin{lem}\label{lemma}
	Let $f \in X_\B$.
\begin{enumerate}
	\item If $f$ is slice preserving, then
\[\|f\|_{X_{\B_I}}=\|f\|_{X_{\B_J}} \quad \text{ for any $I,J \in \s$}.\] 
	\item In the general case, there exists a constant $c$ which depends only on the norm of the operator $\mathcal J$, such that for all $f\in X_\B$
	\[ \|f\|_{X_{\B_I}}\le \|f\|_{X_\B}\le c \|f\|_{X_{\B_I}} \quad \text{for any $I\in \s$.}\]
\end{enumerate}

\end{lem}	  
\begin{proof}
Let us first prove part $(1)$. Let $f$ be slice preserving. Then, the splitting components of $f$ are the same on each slice as holomorphic functions on $\D$: for any $I\in \s$ and any $J\in \s$ orthogonal to $I$, $F=f_I$ and $G\equiv 0$. Hence, 
\[\|f\|_{X_{\B_I}}=\|F\|_{X_\D} \quad \text{for any $I$}.\]

Let us now prove part $(2)$. The first inequality is obvious.\\
For the second one, consider the decomposition (see, e.g., \cite[Lemma 6.11]{GMP})
\[f=f_0+f_1I+f_2J+f_3IJ\]
where $I$ is any imaginary unit, $J$ is any imaginary unit orthogonal to $I$, and $f_0,f_1,f_2,f_3$ are slice preserving functions. By the triangle inequality 
\begin{equation}\label{dis1}
	\|f\|_{X_\B}\le \|f_0\|_{X_\B}+\|f_1\|_{X_\B}+\|f_2\|_{X_\B}+\|f_3\|_{X_\B}=\|f_0\|_{X_{\B_I}}+\|f_1\|_{X_{\B_I}}+\|f_2\|_{X_{\B_I}}+\|f_3\|_{X_{\B_I}}
	\end{equation}
where the last equality follows from part $(1)$.
Via the natural isomorphism $x+yI\mapsto x+yi$, we can identify each $f_\ell$ with a holomorphic function from $\D$ to $\C$ (which preserves the real axis), and the same can be done for the functions $F:=f_0+f_1I$ and $G:=f_2+f_3I$, which are the splitting components of $f$ with respect to the imaginary units $I,J$. 
Moreover, using the operator $\mathcal J$, we can express $f_0,f_1$ as
\[f_0=\frac{F+\mathcal J(F)}{2}, \quad f_1= \frac{F-\mathcal J(F)}{2I}\]
and $f_2,f_3$ as
\[f_2=\frac{G+\mathcal J(G)}{2}, \quad f_3= \frac{G-\mathcal J(G)}{2I}.\]
Hence
\begin{equation}\label{dis2}
	\|f_0\|_{X_{\B_I}}\le \frac{1+\|\mathcal J\|_{\mathcal B(X_\D)}}{2}\|F\|_{X_\D}, \quad \|f_1\|_{X_{\B_I}}\le \frac{1+\|\mathcal J\|_{\mathcal B(X_\D)}}{2}\|F\|_{X_\D} \end{equation}
and
\begin{equation}\label{dis3}
	\|f_2\|_{X_{\B_I}}\le \frac{1+\|\mathcal J\|_{\mathcal B(X_\D)}}{2}\|G\|_{X_\D}, \quad \|f_3\|_{X_{\B_I}}\le \frac{1+\|\mathcal J\|_{\mathcal B(X_\D)}}{2}\|G\|_{X_\D},
	\end{equation}
where the operator norm $\|\mathcal J\|_{\mathcal B(X_\D)}$ is finite by hypothesis.
Moreover, since $f=F+GJ$ on $\B_I$, we also have
\begin{equation}\label{dis4}
	\| F\|_{X_{\D}}=\| F\|_{X_{\B_I}}\le \| f\|_{X_{\B_I}}, \quad \text{and} \quad \| G\|_{X_{\D}}=\| G\|_{X_{\B_I}}\le \| f\|_{X_{\B_I}}.
	\end{equation}
Combining inequalities \eqref{dis1},\eqref{dis2},\eqref{dis3}, and \eqref{dis4}, we conclude that
\[ \|f\|_{X_\B}\le {2}(1+\|\mathcal J\|_{\mathcal B(X_\D)}) \|f\|_{X_{\B_I}}.\]	
\end{proof}

\section{Carleson measures}\label{sec3}
\noindent Let us recall the definition of Carleson measures. 
\begin{defn}
	A finite nonnegative Borel measure $\mu$ on $\B$ is called a {\em $p$-Carleson measure} for the quaternionic Banach space $X_\B$ if it satisfies the embedding inequality
	\begin{equation}\label{mufin}
		\left(\int_{\B}|f(q)|^pd\mu(q)\right)^{1/p}\le c(\mu)\| f\|_{X_\B}
		\end{equation}
	for every $f\in X_\B$, with a constant $c(\mu)$ depending on $\mu$ alone.
	
\end{defn}

\begin{oss}
Notice that whenever the space $X_\B$ contains constant functions, if a measure satisfies the inequality \eqref{mufin} it is necessarily finite. 
\end{oss}

In order to work with the slice structure of $\B$ and $X_\B$, we suitably decompose measures on the quaternionic unit ball (see \cite{Hankel}). First, we can always write a measure $\mu$ on the unit ball $\B$ as $\mu=\mu_\rr+\tilde\mu$ where {$\mu_\rr$ and $\tilde{\mu}$ are nonnegative Borel measures such that} $\tilde \mu(\B \cap \rr)=0$ and $\supp(\mu_\rr)\subseteq \B \cap \rr$. It is non difficult to show that $\mu$ is a 
$p$-Carleson measure if and only if the measures $\mu_\rr$ and $\tilde \mu$ are $p$-Carleson as well.

A further step is guaranteed by the Disintegration Theorem (see, e.g., Theorem 2.28 in \cite{ambrosio}). In fact, in view of this result, any finite measure $\mu$ on $\B$, such that $\mu(\B\cap\rr)=0$,  
can be uniquely decomposed as
\begin{equation}\label{radon}
	d\mu(x+yI)=
	d\mu^+_I(x+yI)d\nu(I)
\end{equation}
where $\nu$ is the measure on the sphere $\s$ defined by 
\[\nu(E)=\mu\left(\{x+yI \in \B\ |\ y>0 \text{ and } I\in E\}\right)\] 
\noindent and $\mu^+_I$ is a (probability) measure on $\B^+_I=\{x+yI\in \B_I \, : \, y \ge 0\}$. 
Hence we can write 
\[\int_{\B}\varphi(x+yI)d\mu(x+yI) = \int_{\s}\int_{\B^+_{I}}\varphi(x+yI) d\mu^+_{I}(x+yI) d\nu(I)\]
for any continuous $\varphi: \B \to \HH$. If in general $\mu$ is a finite measure on $\B$, 
we can decompose $\mu=\mu_\rr+\tilde \mu$, where $\tilde \mu(\B\cap \rr)=0$ and $\supp \mu_\rr \subseteq \B\cap \rr$, 
so that
\begin{equation*}\label{krypton}\int_{\B}\varphi(x+yI)d\mu(x+yI) = \int_{\B\cap \rr}\varphi(x)d\mu_{\rr}(x)+\int_{\s}\int_{\B^+_{I}}\varphi(x+yI) d\tilde\mu^+_I(x+yI)d\nu(I),\end{equation*}
where $\tilde\mu^+_I$ is obtained from $\tilde\mu$ using the Disintegration Theorem.

\subsection{Main result}
Given a nonnegative measure $\mu$ on $\B$,
we can define a nonnegative measure $\mu^s$ on the complex unit disc, by "projecting" all the mass of $\mu$ on a fixed slice. 
 More precisely, suppose first that $\mu(\B\cap \rr)=0$. 
\noindent If $E \subseteq \D$, we define
\[\mu^s(E)=\mu(\{x+yI : \ y>0, \ I\in \s, \  x+yi \in E\}),\]
or, equivalently, 	for any $F$ measurable function defined on $\D\simeq \B_i$, then
\begin{equation}\label{pushmu}
	\int_{\D}F(x+yi)d{\mu^s}(x+yi)=\int_{{\B}}F(x+yi)d\mu(x+yI),
	\end{equation}
i.e. $\mu^s$ is the pushforward of the measure $\mu$ via the projection map $T:\mathbb B\setminus \mathbb R \to \D$, defined by $T(x+yI)=x+yi$, for every $x+yI \in \mathbb B$ such that $y>0$.

Fix now
the imaginary unit $i$
and define ${\mu^+}^{proj}_I$ to be the projection on the half-disc $\B^+_{i}$ of the measure $\mu^+_I$ defined in Equation \eqref{radon}:
\[{\mu^+}_I^{proj}(E)=\mu^+_I(\{x+yI \, \colon \, y>0 \text{ and } x+y i \in E\}), \ \ d{\mu^+}^{proj}_I(x+yi)=d\mu^+_I(x+yI)\] 
for any $E\subseteq \B^+_{i}$. Again, ${\mu^+}^{proj}_I$ is the pushforward of the measure $\mu^+_I$ via the projection $x+yI \mapsto x+yi$ form $\mathbb B_I^+$ to $\B^+_{i}$.
Then, decomposing $\mu$ as in Equation \eqref{radon},
for any $F$ measurable function on $\D^+\simeq \B^+_i$, we have
\begin{equation}\label{murep}
	\begin{aligned}
	&\int_{\B^+_{i}}F(x+yi)d{\mu^s}(x+yi)=\int_{\B}F(T(x+yI))d\mu(x+yI)\\
	&=\int_{\s}\int_{\B^+_{I}}F(T(x+yI)) d\mu^+_{I}(x+yI) d\nu(I)=\int_{\s}\int_{\B^+_{i}}F(x+yi) d{\mu^+}^{proj}_{I}(x+yi) d\nu(I)\\&=\int_{\B^+_{i}}F(x+yi)\int_{\s}d{\mu^+}^{proj}_I(x+yi)d\nu(I).
	\end{aligned}
\end{equation}

\noindent If $\supp(\mu)\cap \rr\neq \emptyset$, we split $\mu$ as $\mu=\mu_\rr+\tilde{\mu}$, with $\supp(\mu_\rr)\subset \rr$ and $\supp(\tilde{\mu})\cap \rr = \emptyset$. Then, we define 
\[\mu^s=\mu_\rr+\tilde{\mu}^s.\] 

\begin{teo}\label{main}
	Let $p\ge 1$.
	If the operator $\mathcal J$ is bounded on the complex Banach space $X_\D$, then
	a finite positive Borel measure $\mu$ is $p$-Carleson for the quaternionic lift $X_\B$, if and only if $\mu^s$ is $p$-Carleson for $X_\D$,
	i.e. 
	\[X_\B \subseteq L^p(\B, d\mu) \quad \iff \quad  X_\D \subseteq L^p(\D, d\mu^s).\]

\end{teo}	
\begin{proof}
	Suppose $\mu$ is $p$-Carleson for $X_\B$, and decompose $\mu$ as $\mu=\mu_\rr+\tilde{\mu}$ with $\supp(\mu_\rr)\subset \rr$ and $\supp (\tilde{\mu})\cap \rr= \emptyset$.
	
	Consider first $F\in X_\D$ such that $F(\bar z)=\overline{F(z)}$, that is $F(z)=\sum_{n\ge 0}z^nr_n$ with real coefficients $r_n$. Then, thanks to Lemma \ref{lemma}, part $(1)$, the function $f(q)=\sum_{n\ge 0}q^nr_n$ belongs to $X_\B$, $\|f\|_{X_\B}=\|F\|_{X_\D}$,  and the modulus $|f|$ is constant on spheres of the form $x+y\s$. 
	
	Therefore, recalling Equation \eqref{murep},
	\begin{equation}\label{normap}
		\begin{aligned}
		\int_{\B}|f(q)|^pd\mu(q)&=\int_{\B\cap \rr}|f(x)|^pd\mu_\rr(x)+ \int_{\s}d\nu(I)\int_{\B^+_{I}}|f(x+yI)|^pd \tilde{\mu}^+_I(x+yI)\\
		&=\int_{\B\cap \rr}|f(x)|^pd\mu_\rr(x)+\int_{\B^+_{i}}|f(x+yi)|^p\int_{\s}d\nu(I){d{{\tilde{\mu}^+}_I}}^{proj}(x+yi)\\
		&=\int_{\D\cap \rr}|F(x)|^pd\mu_\rr(x)+\int_{\D^+}|F(z)|^pd\tilde{\mu}^s(z)=\int_{\D}|F(z)|^pd\mu^s,
	\end{aligned}
	\end{equation}
	where we identified $\B^+_{i}$ with $\D^+$.
Since $\mu$ is $p$-Carleson, there exists a constant $c(\mu)$, depending only on $\mu$, such that
	\begin{equation}
		\begin{aligned}
			&	\int_{\D}|F(z)|^pd\mu^s(z)=\int_{\B}|f(q)|^pd\mu(q) \le  c(\mu)^p \|f\|_{X_\B}^p= c(\mu)^p\|F\|_{X_\D}^p.
		\end{aligned}
	\end{equation}
	
		Consider now a generic $F(z)=\sum_{n\ge 0}z^na_n \in X_\D$. If $a_n=r_n+is_n$  with $r_n,s_n\in \rr$ for every $n$, then we can write
	\[F(z)=\sum_{n \ge 0}z^nr_n+i\sum_{n\ge 0}z^ns_n=:F_1(z)+iF_2(z)\]
	with $F_1(\bar z)=\overline{F_1(z)}$ and $F_2(\bar z)=\overline{F_2(z)}$.
	Since $F_1=\frac{F+\mathcal J(F)}{2}$ and $F_2=\frac{F-\mathcal J(F)}{2i}$ and the operator $\mathcal J$ is bounded, we get that $F_1$ and $F_2 \in X_\D$ and 
	\begin{equation}\label{normacomp}
		\|F_\ell\|_{ X_\D}\le\frac{1+\|\mathcal J\|_{X_\D}}{2}\|F\|_{X_\D} \quad \text{for $\ell=1,2$}.
		\end{equation}

By the triangle inequality, and by what we just proved for functions in $X_\D$ with real coefficients, since $p\ge 1$  
we have
\[	\begin{aligned}
	&\int_\D|F(z)|^pd\mu^s(z)\le
	2^{p-1}\left(\int_\D|F_1(z)|^pd\mu^s(z)+\int_\D|F_2(z)|^pd\mu^s(z)\right)
	\\
	&\le 2^{p-1}c(\mu)^p\left(\|F_1\|_{ X_\D}^p+\|F_2\|_{X_\D}^p\right)\le c(\mu)^p(1+\|\mathcal J\|_{X_\D})^p\|F\|^p_{X_\D}
\end{aligned}
\]	

\noindent and thus we can conclude that $\mu^s$ is a $p$-Carleson measure for $ X_\D$.
	
In the other direction, suppose $\mu^s$ is a $p$-Carleson measure for $X_\D$ and first assume that $f\in X_\B$ is slice preserving.
 Then, identifying $\D$ with $\B_i$, and $f:\B_i \to \B_i$ with a holomorphic function $F$ in $X_\D$, 
	we have
		\begin{equation}
				\begin{aligned}
						&	\|f\|^p_{X_\B}= \|F\|^p_{X_\D}\ge \frac{1}{c(\mu^s)^p}\int_{\D}|F|^pd\mu^s=
						\frac{1}{c(\mu^s)^p}\int_\B|f|^pd\mu,
					\end{aligned}
			\end{equation}
			
\noindent 	{ If $f$ is not slice preserving, consider $I,J \in \s$ with $I \perp J$, and let $f=F+GJ$ be the splitting of $f$ on $\B_I$ with respect to $J$. 	
		Recalling inequality \eqref{normacomp}, and using the same notation, we have that 
		\[2\|F\|_{X_\D}\ge \frac{2}{1+\|\mathcal J\|_{X_\D}}(\|F_1\|_{X_\D}+\|F_2\|_{X_\D})\] and thus 
		\[\|F\|^2_{X_\D}\ge \frac{1}{(1+\|\mathcal J\|_{X_\D})^2}(\|F_1\|^2_{X_\D}+\|F_2\|^2_{X_\D}).\]
		Analogously for the function $G$ we have
		\[\|G\|^2_{X_\D}\ge \frac{1}{(1+\|\mathcal J\|_{X_\D})^2}(\|G_1\|^2_{X_\D}+\|G_2\|^2_{X_\D}).\]
		The functions $F_1,F_2,G_1, G_2$ have real coefficients, hence, in view of Theorem \ref{RF}, we can extend them to slice preserving regular functions and decompose $f$ as
		\[f=F_1+F_2I+G_1J+G_2IJ.\]
		Therefore 
	\begin{equation*}
		\begin{aligned}
			&	\|f\|^p_{X_\B}\ge \|f\|^p_{X_{\B_I}}=(\|F\|^2_{X_\D}+\|G\|^2_{X_\D})^{p/2} \\
			&\ge  \left(\frac{1}{(1+\|\mathcal J\|_{X_\D})^2}\right)^{p/2}(\|F_1\|^2_{X_\D}+\|F_2\|^2_{X_\D}+\|G_1\|^2_{X_\D}+\|G_2\|^2_{X_\D})^{p/2}\\
			&\ge\left(\frac{1}{c(\mu^s)^2(1+\|\mathcal J\|_{X_\D})^2}\right)^{p/2}(\|F_1\|^2_{L^p(\B,d\mu)}+\|F_2\|^2_{L^p(\B,d\mu)}+\|G_1\|^2_{L^p(\B,d\mu)}+\|G_2\|^2_{L^p(\B,d\mu)})^{p/2}
					\end{aligned}
	\end{equation*}
	where the last inequality follows from what we just proved for slice preserving functions.
If $p\ge2$, 
	\begin{equation*}
		\begin{aligned}
			&	\|f\|^p_{X_\B}\ge\frac{1}{c(\mu^s)^p(1+\|\mathcal J\|_{X_\D})^p}(\|F_1\|^p_{L^p(\B,d\mu)}+\|F_2\|^p_{L^p(\B,d\mu)}+\|G_1\|^p_{L^p(\B,d\mu)}+\|G_2\|^p_{L^p(\B,d\mu)})\\
			&\ge\frac{4^{{1}-p}}{c(\mu^s)^p(1+\|\mathcal J\|_{X_\D})^p}\|f\|^p_{L^p(\B,d\mu)},
\end{aligned}
\end{equation*}
that is $\mu$ is $p$-Carleson for $X_\B$. 

If $1\le p<2$,

\begin{equation*}
	\begin{aligned}
		&	\|f\|^p_{X_\B}\ge\frac{2^{p-2}}{c(\mu^s)^p(1+\|\mathcal J\|_{X_\D})^p}(\|F_1\|^p_{L^p(\B,d\mu)}+\|F_2\|^p_{L^p(\B,d\mu)}+\|G_1\|^p_{L^p(\B,d\mu)}+\|G_2\|^p_{L^p(\B,d\mu)})\\
		&\ge\frac{2^{p-2}4^{1-p}}{c(\mu^s)^p(1+\|\mathcal J\|_{X_\D})^p}\|f\|^p_{L^p(\B,d\mu)},
			\end{aligned}
\end{equation*}
thus showing that also in this case $\mu$ is $p$-Carleson for $X_{\B}$.}

\end{proof}

\section{Vanishing Carleson measures}\label{sec4}
Let us recall the definition of vanishing Carleson measures.
\begin{defn}
	A finite nonnegative Borel measure $\mu$ on $\B$ is called a {\em vanishing $p$-Carleson measure} for the quaternionic Banach space $X_\B$ if the embedding 
	\[Id: X_{\B} \to L^p(\B,d\mu)\]
	is compact, namely if for any sequence $\{f_n\}_n \subset X_{\B}$ such that $\sup_{n} \|f_n\|_{X_\B}<+\infty$ there exists a subsequence $\{f_{n_k}\}_k$  converging in norm in $L^p(\B, d\mu)$.
\end{defn}

{
	\begin{defn}Let $\nu$ be a nonnegative Borel measure on $\D$ and
define the measure $\hat{\nu}$ as the reflection of the measure $\nu$ with respect to the real axis, 
	\[\hat{\nu}(E)=\nu(\bar{E}),\]
	where $\bar E =\{x-yi \ : \ x+yi \in E\}$.
\end{defn}
}

{\begin{lem}\label{hatmu}
Let $X_\D$ be a complex Banach space {such that the operator $\mathcal J:X_\D \to X_\D$ is bounded} and let $\nu$ be
a nonnegative Borel measure on $\D$. Then $\nu$ is vanishing $p$-Carleson for $X_\D$ if and only if the measure $\hat{\nu}$ is vanishing $p$-Carleson for $X_\D$.
	 \end{lem}
\begin{proof}
	Since $\hat{\cdot}$ is an involution, it suffices to show only one implication. Consider any measurable function $G:\D \to \C$. Then 
	\[
		\begin{aligned}
			\|G\|^p_{L^p(\D,d\hat{\nu})}&=\int_{\D}|G(z)|^pd\hat{\nu}(z)=\int_{\D}|{{G(z)}}|^pd{{\nu}}(\bar{z})=\int_{\D}\left|{G(\bar{z})}\right|^pd{{\nu}}(z)\\
			&=\int_{\D}\left|\overline{G(\bar{z})}\right|^pd{{\nu}}(z)=\int_{\D}\left|\mathcal J(G)({z})\right|^pd{\nu}(z)=	\|\mathcal J(G)\|^p_{L^p(\D,d{\nu})}.
		\end{aligned}.	\]
		Consider now a sequence $F_n \in X_\D$ such that $\|F_n\|_{X_\D} \le C$. Since $\mathcal J$ is bounded on $X_\D$, the sequence $\mathcal J(F_n)$ is uniformly bounded as well, $\|\mathcal J (F_n)\| \le C \|\mathcal J\|_{\mathcal B}$. 
		Since $\nu$ is a vanishing $p$-Carleson measure, after passing to a subsequence, and using the fact that $Id \circ\mathcal J : X_\D\to L^p(\D,d\nu)$ is continuous, we have that both $F_{n}$ and $\mathcal J(F_n)$ converge in norm in $L^p(\D,d\nu)$ to $F$ and $\mathcal J(F)$ respectively (both belonging to $X_\D$.). Therefore 
				\[\|F_n-F\|_{L^p(\D,d\hat\nu)}=\|\mathcal J(F_n)-\mathcal J(F)\|_{L^p(\D,d\nu)}\]
				converges to $0$ as $n$ goes to $+\infty$, that is 
 $\hat{\nu}$ is vanishing $p$-Carleson as well. 
\end{proof}}

\begin{teo}\label{main2}
		Let $p\ge 1$.
	If the operator $\mathcal J$ is bounded on the complex Banach space $X_\D$, then
	a finite positive Borel measure $\mu$ is vanisihing $p$-Carleson for the quaternionic Banach space $X_\B$, if and only if $\mu^s$ is vanishing $p$-Carleson for $X_\D$,
	i.e. 
	\[X_\B \subseteq L^p(\B, d\mu) \quad \text{compactly} \quad \iff \quad  X_\D \subseteq L^p(\D, d\mu^s) \quad \text{compactly} .\]
	
	\end{teo}
	\begin{proof}
Suppose first that $\mu^s$ is vanishing $p$-Carleson for $X_\D$ and consider a sequence $f_n\in X_{\B}$ such that $\|f_n\|_{X_\B}\le c$. Then we can split $f_n$ on the slice $\B_{ i}$ with respect to the orthogonal imaginary unit $j$ as 
\[f_n=F^{ i}_n+G^{ i}_nj \quad \text{with} \quad \|F_n^{i}\|_{X_\D},\|G_n^{i}\|_{X_\D}\le c.\]
Since $\mu^s$ is vanishing $p$-Carleson, up to subsequences, the holomorphic functions $F_n^{i}$ and $G_n^{i}$ converge in $L^p(\B_{i},d\mu^s)$ to holomorphic (since $X_\D$ embeds continuously in $\Hol(\D)) $ functions $F^{i}$ and $G^{i}$. Moreover, thanks to the Representation Formula \ref{RF}, 
the splitting components of $f_n$ on a different slice $\B_I$, with respect to $L\in \s$, orthogonal to $I$, satisfy the following equality:
\[(F_n^I+G_n^IL)(x+yI)=\frac{1-I{i}}{2}(F_n^{i}+G_n^{i}j)(x+y{i})+\frac{1+I{i}}{2}(F_n^{i}+G_n^{i}j)(x-y{i}),\]
 so that the same relation holds for the limit functions
\[(F^I+G^IL)(x+yI)=\frac{1-I{i}}{2}(F^{i}+G^{i}j)(x+y{i})+\frac{1+I{i}}{2}(F^{i}+G^{i}j)(x-y{i})\]  which implies that, up to subsequences, $f_n$ converges uniformly on compact sets to a slice regular function $f$.
We need to verify that \[\|f_n-f\|_{L^p(\B,d\mu)}\to 0 \quad \text{as $n\to \infty$}.\]
Thanks to the Representation Formula and to Equation \eqref{pushmu} we have that
\begin{equation}
	\begin{aligned}
	\|f_n-f\|^p_{L^p(\B,d\mu)}&=\int_{\B}|f_n(x+yI)-f(x+yI)|^pd\mu(x+yI)\\
	&\le c_1(p)\left(\int_{\B}|(f_n-f)(x+yi)|^pd\mu(x+yI)+\int_{\B}|(f_n-f)(x-yi)|^pd\mu(x+yI)\right)\\
	&= c_1(p)\left(\int_{\B_{i}}|(f_n-f)(x+yi)|^pd\mu^s(x+yi)+\int_{\B_{i}}|{(f_n-f)(x-yi)}|^pd\mu^s(x+yi)\right)\\
	&= c_1(p)\int_{\B_{i}}|(f_n-f)(x+yi)|^p(d\mu^s(x+yi)+d\hat{\mu}^s(x+yi))
	\end{aligned}
	\end{equation}
	where $c_1(p)$ is a constant depending only on $p$. 
	Applying the Splitting Lemma, there exists a constant $c_2(p)$, depending only on $p$, such that	
\begin{equation}
	\begin{aligned}
	&	\|f_n-f\|^p_{L^p(\B,d\mu)}\\
		&\le c_1(p)\int_{\B_{i}}\left(|({F_n^{i}}-F^{i})(x+yi)|^2+|({G_n^{i}}-G^{i})(x+yi)|^2\right)^{p/2}
		(d\mu^s+d\hat{\mu}^s)(x+yi)\\
		&\le c_2(p)\int_{\B_{i}}\left(|(F_n^{i}-F^{i})(x+yi)|^p+|(G_n^{i}-G^{i})(x+yi)|^p\right)(d\mu^s+d\hat{\mu}^s)(x+yi)\\
		&=c_2(p)\left(\| F_n^{i}-F^{i}\|^p_{L^p(\B_{i},d\mu^s+d\hat{\mu}^s )}+\| G_n^{i}-G^{i}\|^p_{L^p(\B_{i},d\mu^s+d\hat{\mu}^s )}\right)
	\end{aligned}
\end{equation}
which, after passing to a subsequence, converges to $0$ since, by Lemma \ref{hatmu}, $d\mu^s+d\hat{\mu}^s$ is a vanishing $p$-Carleson measure.

In the other direction, we proceed as in the first part of the proof of Theorem \ref{main}. Suppose that $\mu$ is vanishing $p$-Carleson and consider a uniformly bounded sequence $F_n$ in $X_{\D}$. 
Let $f_n:=\ext(F_n)$ for any $n$. Recalling Lemma \ref{lemma} we have that there exists a constant $c$ such that
\[\|f_n\|_{X_\B}\le c\|F_n\| \quad \text{for every $n$},\] so that $f_n$ is a uniformly bounded sequence in $X_\B$. Since $\mu$ is vanishing $p$-Carleson, and since $X_\B$ embeds continuously in $\mathcal{ SR}(\B)$, up to taking a subsequence, $f_n$ converges to a slice regular function $f$ with respect to the $L^p(\B,d\mu)$-norm. Let $F\in \Hol(\D)$ be such that $f=\ext(F)$. We want to show that $F_n$ converges in $L^p(\D,d\mu^s)$ to the function $F$. 
 
Consider first the case in which $F_n$ are all such that $F_n(\bar z)=\overline{F_n(z)}$.
	Then, recalling the Representation and Extension Formula \ref{RF} and Lemma \ref{lemma}, the functions $f_n(q)=\ext(F_n)(q)$ are slice preserving functions in $X_\B$, with $\|f_n\|_{X_\B}=\|F_n\|_{X_\D}$ and the limit function $f=\ext(F)$ is slice preserving as well. 
Reasoning as in \eqref{normap}, we have that 
\begin{equation}\label{normap2}
	\begin{aligned}
		\int_{\B}|f_n(q)-f(q)|^pd\mu(q)
=\int_{\D}|F_n(z)-F(z)|^pd\mu^s(z).
	\end{aligned}
\end{equation}
Hence, since $\mu$ is vanishing $p$-Carleson, up to taking a subsequence, 
\begin{equation}
	\begin{aligned}
		&\|F_n -F\|_{L^p(\D,d\mu^s)}=\|f_n -f\|_{L^p(\B,d\mu)}	\end{aligned}
\end{equation}
goes to $0$ as $n$ goes to $+\infty$. 

Consider now any $F_n(z)=\sum_{k\ge 0}z^ka^{(n)}_k \in X_\D$. If $a^{(n)}_k=r^{(n)}_k+is^{(n)}_k$  with $r^{(n)}_k,s^{(n)}_k\in \rr$ for every $k$, and every $n$, then we can write
\[F_n(z)=\sum_{k \ge 0}z^kr_k^{(n)}+i\sum_{k\ge 0}z^ks_k^{(n)}=:F_n^1(z)+iF_n^2(z)\]
with $F^1_n(\bar z)=\overline{F_n^1(z)}$ and $F_n^2(\bar z)=\overline{F_n^2(z)}$.
Since $F_n^1=\frac{F_n+\mathcal J(F_n)}{2}$ and $F_n^2=\frac{F_n-\mathcal J(F_n)}{2i}$ and the operator $\mathcal J$ is bounded, we get that $F_n^1$ and $F_n^2 \in X_\D$ for any $n$ and 
\begin{equation}
	\|F_n^\ell\|_{ X_\D}\le\frac{1+\|\mathcal J\|_{X_\D}}{2}\|F_n\|_{X_\D} \quad \text{for $\ell=1,2$},
\end{equation}
so that the sequences $F_n^1$ and $F_n^2$ are both uniformly bounded. 

Let $F^1$ and $F^2$ be the holomorphic functions with real coefficients, such that the limit function $f=\ext(F)$ can be written as $f=\ext(F^1)+i\ext(F^2)$. 
By the triangle inequality, since $p\ge 1$, 
we have then

\[	\begin{aligned}
	&\int_\D|F_n(z)-F|^pd\mu^s(z)\le
	2^{p-1}\left(\int_\D|F_n^1(z)-F^1|^pd\mu^s(z)+\int_\D|F_n^2(z)-F^2(z)|^pd\mu^s(z)\right)
	\\
	&= 2^{p-1}\left(\|F_n^1-F^1\|_{L^p(\D,d\mu^s)}^p+\|F_n^2-F^2\|_{L^p(\D,d\mu^s)}^p\right)
\end{aligned}
\]	
which, up to subsequences, goes to $0$ as $n$ goes to $+\infty$, thanks to what we just proved for sequences of functions in $X_\D$ with real coefficients.
\noindent Thus we can conclude that $\mu^s$ is a vanishing $p$-Carleson measure for $ X_\D$.

	\end{proof}

\end{document}